\documentclass[11pt,a4paper,twoside]{amsart}

\usepackage{graphics}
\usepackage{epic}
\usepackage{amssymb}
\usepackage{amsmath}

\newtheorem{theo}{\bf Theorem}[section]

\newtheorem{lemma}[theo]{\bf Lemma}
\newtheorem{defi}[theo]{\bf Definition}
\newtheorem{coro}[theo]{\bf Corollary}

\newcommand{\dash}{{
    \dashline[50]{2}(0,0)(30,0)
    \dashline[50]{2}(0,-10)(30,-10)
    \dashline[50]{2}(0,-20)(10,-20)
    \dashline[50]{2}(0,0)(0,-20)
    \dashline[50]{2}(10,0)(10,-20)
    \dashline[50]{2}(20,0)(20,-10)
    \dashline[50]{2}(30,0)(30,-10)
  }}

\newcommand{\pair}[2]{
  {\begin{smallmatrix}
      #1 \\
      #2
    \end{smallmatrix}}}

\newcommand{\sarray}[1]{
  {\begin{smallmatrix}
      #1
    \end{smallmatrix}}}

\newcommand{\cupdot}{{\uplus}}
\newcommand{\PS}{{\rm PS}}
\newcommand{\PT}{{\rm PT}}
\newcommand{\ST}{{\rm ST}}
\newcommand{\LHS}{{\rm LHS}}
\newcommand{\RHS}{{\rm RHS}}

\newcommand{\vd}{{\,\vdash\,}}
\newcommand{\sh}[1]{{\,{\rm sh}\,#1}}

\newcommand{\rsgn}[1]{{\,{\rm rsgn}\,#1}}
\newcommand{\invsgn}[1]{{\,{\rm invsgn}\,#1}}
\newcommand{\sgn}[1]{{\,{\rm sgn}\,#1}}

\newcommand{\bigexample}{{
    \begin{figure}
      \fbox{\begin{minipage}{5.28in}
          Let $\alpha=(3,1)$ and $n=3$. There are 10 skew shapes
          $\lambda/\alpha\vd3$. Here we have evaluated
          $(-1)^{v(\lambda)}I_{\lambda/\alpha}^2$ for each one of them:
          \begin{center}
            \setlength{\unitlength}{0.3mm}
            \begin{tabular}{c}
              $-1$ \\
              \begin{picture}(70,60)(-10,-55)
                \dash
                \thicklines
                
                \put(30,0){\line(1,0){30}}
                \put(30,-10){\line(1,0){30}}
                
                \put(30,0){\line(0,-1){10}}
                \put(40,0){\line(0,-1){10}}
                \put(50,0){\line(0,-1){10}}
                \put(60,0){\line(0,-1){10}}
              \end{picture}
            \end{tabular}
            \begin{tabular}{c}
              $+1$ \\
              \begin{picture}(70,60)(-15,-55)
                \dash
                \thicklines
                
                \put(30,0){\line(1,0){20}}
                \put(30,-10){\line(1,0){20}}
                \put(10,-10){\line(1,0){10}}
                \put(10,-20){\line(1,0){10}}
                
                \put(30,0){\line(0,-1){10}}
                \put(40,0){\line(0,-1){10}}
                \put(50,0){\line(0,-1){10}}
                \put(10,-10){\line(0,-1){10}}
                \put(20,-10){\line(0,-1){10}}
              \end{picture}
            \end{tabular}
            \begin{tabular}{c}
              $-1$ \\
              \begin{picture}(70,60)(-15,-55)
                \dash
                \thicklines
                
                \put(30,0){\line(1,0){20}}
                \put(30,-10){\line(1,0){20}}
                \put(0,-20){\line(1,0){10}}
                \put(0,-30){\line(1,0){10}}

                \put(30,0){\line(0,-1){10}}
                \put(40,0){\line(0,-1){10}}
                \put(50,0){\line(0,-1){10}}
                \put(0,-20){\line(0,-1){10}}
                \put(10,-20){\line(0,-1){10}}
              \end{picture}
            \end{tabular}
            \begin{tabular}{c}
              $-1$ \\
              \begin{picture}(70,60)(-20,-55)
                \dash
                \thicklines

                \put(30,0){\line(1,0){10}}
                \put(10,-10){\line(1,0){30}}
                \put(10,-20){\line(1,0){20}}

                \put(10,-10){\line(0,-1){10}}
                \put(20,-10){\line(0,-1){10}}
                \put(30,0){\line(0,-1){20}}
                \put(40,0){\line(0,-1){10}}
              \end{picture}
            \end{tabular}
            \begin{tabular}{c}
              $+0$ \\
              \begin{picture}(70,60)(-20,-55)
                \dash
                \thicklines

                \put(30,0){\line(1,0){10}}
                \put(30,-10){\line(1,0){10}}
                \put(10,-10){\line(1,0){10}}
                \put(0,-20){\line(1,0){20}}
                \put(0,-30){\line(1,0){10}}

                \put(0,-20){\line(0,-1){10}}
                \put(10,-10){\line(0,-1){20}}
                \put(20,-10){\line(0,-1){10}}
                \put(30,0){\line(0,-1){10}}
                \put(40,0){\line(0,-1){10}}
              \end{picture}
            \end{tabular}
            \begin{tabular}{c}
              $+1$ \\
              \begin{picture}(70,60)(-20,-55)
                \dash
                \thicklines

                \put(30,0){\line(1,0){10}}
                \put(30,-10){\line(1,0){10}}
                \put(0,-20){\line(1,0){10}}
                \put(0,-30){\line(1,0){10}}
                \put(0,-40){\line(1,0){10}}

                \put(0,-20){\line(0,-1){20}}
                \put(10,-20){\line(0,-1){20}}
                \put(30,0){\line(0,-1){10}}
                \put(40,0){\line(0,-1){10}}
              \end{picture}
            \end{tabular}
            \begin{tabular}{c}
              $-1$ \\
              \begin{picture}(70,60)(-25,-55)
                \dash
                \thicklines

                \put(10,-10){\line(1,0){20}}
                \put(0,-20){\line(1,0){30}}
                \put(0,-30){\line(1,0){10}}

                \put(0,-20){\line(0,-1){10}}
                \put(10,-10){\line(0,-1){20}}
                \put(20,-10){\line(0,-1){10}}
                \put(30,-10){\line(0,-1){10}}
              \end{picture}
            \end{tabular}
            \begin{tabular}{c}
              $+0$ \\
              \begin{picture}(75,60)(-25,-55)
                \dash
                \thicklines

                \put(10,-10){\line(1,0){10}}
                \put(0,-20){\line(1,0){20}}
                \put(0,-30){\line(1,0){20}}

                \put(0,-20){\line(0,-1){10}}
                \put(10,-10){\line(0,-1){20}}
                \put(20,-10){\line(0,-1){20}}
              \end{picture}
            \end{tabular}
            \begin{tabular}{c}
              $-1$ \\
              \begin{picture}(70,60)(-25,-55)
                \dash
                \thicklines

                \put(10,-10){\line(1,0){10}}
                \put(0,-20){\line(1,0){20}}
                \put(0,-30){\line(1,0){10}}
                \put(0,-40){\line(1,0){10}}

                \put(0,-20){\line(0,-1){20}}
                \put(10,-10){\line(0,-1){30}}
                \put(20,-10){\line(0,-1){10}}
              \end{picture}
            \end{tabular}
            \begin{tabular}{c}
              $+1$ \\
              \begin{picture}(70,60)(-25,-55)
                \dash
                \thicklines

                \put(0,-20){\line(1,0){10}}
                \put(0,-30){\line(1,0){10}}
                \put(0,-40){\line(1,0){10}}
                \put(0,-50){\line(1,0){10}}

                \put(0,-20){\line(0,-1){30}}
                \put(10,-20){\line(0,-1){30}}
              \end{picture}
            \end{tabular}
          \end{center}
          (It happens that all these skew shapes have sign-imbalance 0 or 1, but
          in larger examples we would find much more exotic integers, like $-7$
          for instance.) Now we compute $(-1)^{v(\mu)}I_{\alpha/\mu}^2$
          for the two skew shapes $\alpha/\mu\vd2$:
          \begin{center}
            \begin{tabular}{c}
              $-1$ \\
              \begin{picture}(75,30)(-25,-25)
                \dash
                \thicklines

                \put(10,0){\line(1,0){20}}
                \put(10,-10){\line(1,0){20}}

                \put(10,0){\line(0,-1){10}}
                \put(20,0){\line(0,-1){10}}
                \put(30,0){\line(0,-1){10}}
              \end{picture}
            \end{tabular}
            \begin{tabular}{c}
              $+0$ \\
              \begin{picture}(75,30)(-25,-25)
                \dash
                \thicklines

                \put(20,0){\line(1,0){10}}
                \put(20,-10){\line(1,0){10}}
                \put(0,-10){\line(1,0){10}}
                \put(0,-20){\line(1,0){10}}

                \put(0,-10){\line(0,-1){10}}
                \put(10,-10){\line(0,-1){10}}
                \put(20,0){\line(0,-1){10}}
                \put(30,0){\line(0,-1){10}}
              \end{picture}
            \end{tabular}
          \end{center}
          Finally, there is only one skew shape $\alpha/\mu\vd3$:
          \begin{center}
            \begin{tabular}{c}
              $+1$ \\
              \begin{picture}(75,30)(-25,-25)
                \dash
                \thicklines

                \put(10,0){\line(1,0){20}}
                \put(0,-10){\line(1,0){30}}
                \put(0,-20){\line(1,0){10}}

                \put(0,-10){\line(0,-1){10}}
                \put(10,0){\line(0,-1){20}}
                \put(20,0){\line(0,-1){10}}
                \put(30,0){\line(0,-1){10}}
              \end{picture}
            \end{tabular}
          \end{center}
          We check that
          $$
          \sum_{\lambda/\alpha\vd3}
          (-1)^{v(\lambda)}I_{\lambda/\alpha}^2=-2=-1-1=
          \sum_{\alpha/\mu\vd2}
          (-1)^{v(\mu)}I_{\alpha/\mu}^2
          -\sum_{\alpha/\mu\vd3}
          (-1)^{v(\mu)}I_{\alpha/\mu}^2.
          $$
        \end{minipage}}
      \caption{Example of Theorem~\ref{th:inout}.}
      \label{fig:bigexample}
    \end{figure}
  }}

\begin{document}

\title{On the sign-imbalance of skew partition shapes}
\author{Jonas Sj{\"o}strand}
\address{Department of Mathematics, Royal Institute of Technology \\
  SE-100 44 Stockholm, Sweden}
\email{jonass@kth.se}

\begin{abstract}
  Let the {\em sign} of a skew standard Young tableau be the sign of the
  permutation you get by reading it row by row from left to right,
  like a book. We examine how the sign property is transferred by the
  skew Robinson-Schensted correspondence invented by Sagan and Stanley.
  The result is a remarkably simple generalization of the
  ordinary non-skew formula.

  The sum of the signs of all standard tableaux on
  a given skew shape is the {\em sign-imbalance} of that shape.
  We generalize previous results on the
  sign-imbalance of ordinary partition shapes to
  skew ones.
\end{abstract}


\date{22 July 2005}

\maketitle

\section{Introduction}
\noindent
A {\em labelled poset} $(P,\omega)$ is an $n$-element poset $P$ with
a bijection $\omega\ :\ P\rightarrow[n]=\{1,2,\ldots,n\}$
called the {\em labelling} of $P$. A {\em linear extension}
of $P$ is an order-preserving bijection $f\ :\ P\rightarrow[n]$.
It is natural to define the {\em sign} of $f$ as
$-1$ to the power of the number of inversions with respect to
the labelling, i.e.,~pairs $x,y\in P$ such that
$\omega(x)<\omega(y)$ and $f(x)>f(y)$.
The {\em sign-imbalance} $I_{P,\omega}$ of $(P,\omega)$ is the sum of the
signs of all linear extensions of $P$.
Note that $I_{P,\omega}$ is independent of the labelling $\omega$
up to sign. In this paper we will mainly discuss the square
of sign-imbalances, and then we may drop the $\omega$ and write
$I_P^2=I_{P,\omega}^2$.

If $I_P^2=0$ the poset is {\em sign-balanced}. Such posets have
been studied since 1989 by F.~Ruskey \cite{ruskey2}, \cite{ruskey1},
R.~Stanley \cite{stanley}, and
D.~White \cite{white}.
It is a vast subject however, and most of the work has been
devoted to a certain class of posets: the partition
shapes (or Young diagrams). Though no one so far has been able to
completely characterize the sign-balanced partition shapes,
this research direction has offered a lot of interesting results.
Many people have studied the more general notion of
sign-imbalance of partition shapes, among those
T.~Lam \cite{lam}, A.~Reifegerste \cite{astrid},
J.~Sj\"ostrand \cite{jonas},
M.~Shimozono and D.~White \cite{shimozonowhite}, 
R.~Stanley \cite{stanley}, and 
D.~White \cite{white}.

Young tableaux play a central role in the
theory of symmetric functions (see~\cite{fulton}) and there
are lots of useful tools for working with them that are
not applicable to general posets. One outstanding tool is the
Robinson-Schensted correspondence which has produced nice results
also in the field of sign-imbalance, see \cite{jonas}, \cite{astrid},
and \cite{shimozonowhite}.

As suggested in \cite{jonas} a natural step from partition shapes
towards more general posets would be to study {\em skew}
partition shapes. They have the advantage of being surrounded by
a well-known algebraic and combinatorial machinery
just like the ordinary shapes, and possibly they might shed some light
on the sign-imbalance of the latter ones as well.
We will use
a generalization of the Robinson-Schensted algorithm for skew
tableaux by B.~Sagan and R.~Stanley \cite{saganstanley}.

In a recent paper~\cite[Theorem~4.3 and~5.7]{frank}
E.~Soprunova and F.~Sottile
show that $|I_{P,\omega}|$ is a lower bound for the number
of real solutions to certain polynomial systems. Theorem~6.4
in~\cite{frank} says that $|I_{P,\omega}|$ is
the characteristic of the Wronski projection on certain projective
varieties associated with $P$. When $P$ is a skew partition shape
this is applicable to skew Schubert varieties in Grassmanians
(Richardson varieties).

An outline of this paper:
\begin{itemize}
\item After some basic definitions in section~\ref{sec:prel},
  in section~\ref{sec:skewrs}
  we briefly recall Sagan and Stanley's skew RS-correspondence
  from~\cite{saganstanley}.
\item
  In section~\ref{sec:results}
  we state our main results without proofs
  and examine their connection to old results.
\item
  In section~\ref{sec:mainproof} and ~\ref{sec:inoutproof}
  we prove our main theorems through
  a straightforward but technical analysis.
\item
  In section~\ref{sec:specialisations} we
  examine a couple of interesting corollaries
  to our main results. One corollary is
  a surprising formula for
  the square of the sign-imbalance of any ordinary shape.
\item
  Finally, in section~\ref{sec:future} we
  suggest some future research directions.
\end{itemize}

\section{Preliminaries\label{sec:prel}}
\noindent
An (ordinary) {\em $n$-shape} $\lambda=(\lambda_1,\lambda_2,\ldots)$ is
a graphical representation (a Ferrers diagram)
of an integer partition of $n=\sum_i\lambda_i$. We write
$\lambda\vd n$ or $|\lambda|=n$. The {\em coordinates} of a cell
is the pair $(r,c)$ where $r$ and $c$ are the row and column indices.
Example:
\begin{center}
  \setlength{\unitlength}{0.5mm}
  \begin{picture}(50,40)(0,-45)
    \put(-70,-27){$(6,4,2,2,1)\ \ =$}
    \put(0,0){\line(1,0){60}}
    \put(0,-10){\line(1,0){60}}
    \put(0,-20){\line(1,0){40}}
    \put(0,-30){\line(1,0){20}}
    \put(0,-40){\line(1,0){20}}
    \put(0,-50){\line(1,0){10}}
    \put(0,0){\line(0,-1){50}}
    \put(10,0){\line(0,-1){50}}
    \put(20,0){\line(0,-1){40}}
    \put(30,0){\line(0,-1){20}}
    \put(40,0){\line(0,-1){20}}
    \put(50,0){\line(0,-1){10}}
    \put(60,0){\line(0,-1){10}}
    \put(38,-32){$(3,2)$}
    \put(35,-30){\vector(-4,1){20}}
  \end{picture}
\end{center}

A shape $\mu$ is a {\em subshape} of a shape $\lambda$ if
$\mu_i\leq\lambda_i$ for all $i$.
For any subshape $\mu\subseteq\lambda$ the
{\em skew} shape $\lambda/\mu$ is $\lambda$ with $\mu$ deleted.
A {\em skew $n$-shape} $\lambda/\mu$ is a skew shape with $n$ cells, and
we write $\lambda/\mu\vd n$ or $|\lambda/\mu|=n$.
Here is an example of a skew 6-shape:
\begin{center}
  \setlength{\unitlength}{0.5mm}
  \begin{picture}(50,50)(0,-50)
    \put(-100,-27){$(6,4,2,2,1)/(4,3,2)\ \ =$}

    \dashline[-7]{1}(0,0)(60,0)
    \dashline[-7]{1}(0,-10)(60,-10)
    \dashline[-7]{1}(0,-20)(40,-20)
    \dashline[-7]{1}(0,-30)(20,-30)
    \dashline[-7]{1}(0,-40)(20,-40)
    \dashline[-7]{1}(0,-50)(10,-50)
    \dashline[-7]{1}(0,0)(0,-50)
    \dashline[-7]{1}(10,0)(10,-50)
    \dashline[-7]{1}(20,0)(20,-40)
    \dashline[-7]{1}(30,0)(30,-20)
    \dashline[-7]{1}(40,0)(40,-20)
    \dashline[-7]{1}(50,0)(50,-10)
    \dashline[-7]{1}(60,0)(60,-10)

    \put(40,0){\line(1,0){20}}
    \put(30,-10){\line(1,0){30}}
    \put(30,-20){\line(1,0){10}}
    \put(0,-30){\line(1,0){20}}
    \put(0,-40){\line(1,0){20}}
    \put(0,-50){\line(1,0){10}}
    \put(0,-30){\line(0,-1){20}}
    \put(10,-30){\line(0,-1){20}}
    \put(30,-10){\line(0,-1){10}}
    \put(20,-30){\line(0,-1){10}}
    \put(40,0){\line(0,-1){20}}
    \put(50,0){\line(0,-1){10}}
    \put(60,0){\line(0,-1){10}}
  \end{picture}
\end{center}

A {\em domino} is a rectangle consisting
of two cells. For an ordinary shape $\lambda$, let
$v(\lambda)$ denote the maximal number of disjoint
vertical dominoes that fit in the
shape $\lambda$.

A {\em (partial) tableau} $T$ on a skew $n$-shape
$\lambda/\mu$ is a labelling of the cells
of $\lambda/\mu$ with $n$ distinct real numbers
such that every number is greater
than its neighbours above and to the left.
We let $\sharp T=n$ denote the number of entries
in $T$, and $\PT(\lambda/\mu)$ denote the
set of partial tableaux on $\lambda/\mu$.

A {\em standard tableau} on a skew $n$-shape
is a tableau with the numbers $[n]=\{1,2,\ldots,n\}$.
We let $\ST(\lambda/\mu)$ denote the set of
standard tableaux on the shape $\lambda/\mu$.
Here is an example:
\begin{center}
  \setlength{\unitlength}{0.5mm}
  \begin{picture}(50,50)(0,-50)
    \dashline[-7]{1}(0,0)(60,0)
    \dashline[-7]{1}(0,-10)(60,-10)
    \dashline[-7]{1}(0,-20)(40,-20)
    \dashline[-7]{1}(0,-30)(20,-30)
    \dashline[-7]{1}(0,-40)(20,-40)
    \dashline[-7]{1}(0,-50)(10,-50)
    \dashline[-7]{1}(0,0)(0,-50)
    \dashline[-7]{1}(10,0)(10,-50)
    \dashline[-7]{1}(20,0)(20,-40)
    \dashline[-7]{1}(30,0)(30,-20)
    \dashline[-7]{1}(40,0)(40,-20)
    \dashline[-7]{1}(50,0)(50,-10)
    \dashline[-7]{1}(60,0)(60,-10)

    \put(40,0){\line(1,0){20}}
    \put(30,-10){\line(1,0){30}}
    \put(30,-20){\line(1,0){10}}
    \put(0,-30){\line(1,0){20}}
    \put(0,-40){\line(1,0){20}}
    \put(0,-50){\line(1,0){10}}
    \put(0,-30){\line(0,-1){20}}
    \put(10,-30){\line(0,-1){20}}
    \put(30,-10){\line(0,-1){10}}
    \put(20,-30){\line(0,-1){10}}
    \put(40,0){\line(0,-1){20}}
    \put(50,0){\line(0,-1){10}}
    \put(60,0){\line(0,-1){10}}
    \put(40,-10){\makebox(10,10){1}}
    \put(50,-10){\makebox(10,10){4}}
    \put(30,-20){\makebox(10,10){3}}
    \put(0,-40){\makebox(10,10){2}}
    \put(10,-40){\makebox(10,10){6}}
    \put(0,-50){\makebox(10,10){5}}
  \end{picture}
\end{center}
The (skew) shape of a tableau $T$ is denoted by $\sh T$.
Note that it is not sufficient to look at the cells of $T$
in order to determine its shape; we must think of the tableau
as remembering its underlying skew shape. (For instance,
$(6,4,2,2,1)/(4,3,2)$ and $(6,4,3,2,1)/(4,3,3)$ are distinct
skew shapes that have the same set of cells.)

The {\em sign} of a number sequence $w_1w_2\cdots w_k$
is $(-1)^{\sharp\{(i,j)\,:\,i<j,\,w_i>w_j\}}$, so it is
$+1$ for an even number of inversions, $-1$ otherwise.
The {\em inverse sign} is defined to be
$(-1)^{\sharp\{(i,j)\,:\,i<j,\,w_i<w_j\}}$.

The {\em sign} $\sgn T$ and the {\em inverse sign} $\invsgn T$
of a tableau $T$ are the sign respectively the inverse sign
of the sequence you get
by reading the entries row by row, from left to right and from
top to bottom, like a book.
Our example tableau has 4 inversions and 11 non-inversions,
so $\sgn T=+1$ and $\invsgn T=-1$.
\begin{defi}
  The {\em sign-imbalance} $I_{\lambda/\mu}$ of a skew shape
  $\lambda/\mu$ is the sum of the signs of all standard tableaux
  on that shape:
  $$I_{\lambda/\mu}=\sum_{T\in\ST(\lambda/\mu)}\sgn T.$$
\end{defi}
\noindent
An empty tableau has positive sign and
$I_{\lambda/\lambda}=I_\emptyset=1$.

A {\em biword} $\pi$ is a sequence of vertical pairs of positive
integers $\pi=\pair{i_1i_2\cdots i_k}{j_1j_2\cdots j_k}$
with $i_1\le i_2\le\cdots\le i_k$. We define the top and bottom
lines of $\pi$ by $\hat\pi=i_1i_2\cdots i_k$ and
$\check\pi=j_1j_2\cdots j_k$.
A {\em partial $n$-permutation} is a biword where in each line
the elements are distinct and of size at most $n$. Let
$\PS_n$ denote the set of partial $n$-permutations.

For each $\pi\in\PS_n$ we associate an ordinary
$n$-permutation $\bar{\pi}\in S_n$
constructed as follows: First take the numbers among $1,2,\ldots,n$
that do not belong to $\hat\pi$ and sort them in increasing order
$a_1<a_2<\cdots<a_\ell$.
Then sort the numbers among $1,2,\ldots,n$
that do not belong to $\check\pi$ in increasing order
$b_1<b_2<\cdots<b_\ell$.
Now insert the vertical pairs $\pair{a_r}{b_r}$, $1\le r\le\ell$
into $\pi$ so that the top line remains increasingly ordered
(and hence must be $12\cdots n$). The bottom line is
a permutation (in single-row notation) which we denote $\bar\pi$.
Example: If $n=5$ and $\pi=\pair{124}{423}$ then $\bar\pi=42135$.

In the following we let $\cupdot$ denote disjoint union interpreted
liberally. For instance, we will write $\check\pi\cupdot T=[n]$
meaning that the set of numbers appearing in $\check\pi$ and the
set of entries of the tableau $T$ are disjoint and their union is $[n]$.

\section{The skew RS-correspondence\label{sec:skewrs}}
In~\cite{saganstanley} Bruce Sagan and Richard Stanley
introduced several analogues of the
Robinson-Schensted algorithm for skew Young tableaux. Their main
result is the following theorem.
\begin{theo}[Sagan and Stanley; 1990]\label{th:ss}
  Let $n$ be a fixed positive integer and $\alpha$ a fixed partition
  (not necessarily of $n$). Then there is bijection
  $$(\pi,T,U)\leftrightarrow(P,Q)$$
  between $\pi\in\PS_n$
  with $T,U\in\PT(\alpha/\mu)$ such that
  $\check\pi\cupdot T=\hat\pi\cupdot U=[n]$, on the one hand, and
  $P,Q\in\ST(\lambda/\alpha)$ such that $\lambda/\alpha\vd n$,
  on the other.
\end{theo}
\noindent
Though we will assume detailed familiarity with it,
we do not define the bijection here, but refer to \cite{saganstanley}
for the original presentation.


\section{Our results\label{sec:results}}
\noindent
In~\cite{jonas} and \cite{astrid} the author and Astrid Reifegerste
independently discovered the formula for sign transfer under
the RS-correspondence:
\begin{theo}[Reifegerste; Sj\"ostrand; 2003]\label{th:jonasastrid}
  Under the (ordinary) RS-correspondence $\pi\leftrightarrow(P,Q)$ we have
  $$\sgn\pi=(-1)^{v(\lambda)}\sgn P\sgn Q$$
  where $\lambda$ is the shape of $P$ and $Q$.
\end{theo}
\noindent
Our main theorem is a generalization of this
to Sagan and Stanley's skew RS-correspondence:
\begin{theo}\label{th:main}
  Under the skew RS-correspondence $(\pi,T,U)\leftrightarrow(P,Q)$
  we have
  $$(-1)^{v(\lambda)}\sgn P\sgn Q=(-1)^{|\alpha|}(-1)^{v(\mu)+|\mu|}
  \sgn T\sgn U\sgn{\bar\pi}$$
  where $\sh P=\sh Q=\lambda/\alpha$ and $\sh T=\sh U=\alpha/\mu$.
\end{theo}
\noindent
Note that if
$\alpha=\emptyset$ the theorem reduces to Theorem~\ref{th:jonasastrid}.
\vspace{2mm}

\noindent
{\bf Remark.}
If we specialise to the skew RS-correspondence
$(\pi,T)\leftrightarrow P$ of involutions
(see Corollary~3.4 in~\cite{saganstanley}),
Theorem~\ref{th:main} gives that
$$(-1)^{v(\lambda)}=\sgn{\bar\pi}(-1)^{v(\mu)+|\mu|+|\alpha|},$$
where $\sh P=\lambda/\alpha$ and $\sh T=\alpha/\mu$.
This is also a simple consequence of Corollary~3.6
in~\cite{saganstanley} which is a generalization of a theorem by
Sch\"utzenberger \cite[page~127]{schutzenberger}
(see also \cite[exercise~7.28~a]{enum2}).
\vspace{2mm}

A fundamental application of Theorem~\ref{th:jonasastrid}
appearing in both \cite{jonas} and \cite{astrid} is the following
theorem.
\bigexample
\begin{theo}[Reifegerste; Sj\"ostrand; 2003]\label{th:jonasastrid2}
  For all $n\geq2$
  $$\sum_{\lambda\vdash n}(-1)^{v(\lambda)}I_\lambda^2=0.$$
\end{theo}
\noindent
We give a natural generalization of this using
Theorem~\ref{th:main}. It may be called a
``sign-imbalance analogue'' to Corollary~2.2 in~\cite{saganstanley}.
\begin{theo}\label{th:inout}
  Let $\alpha$ be a fixed partition and let $n$ be a positive integer.
  Then
  $$
  \sum_{\lambda/\alpha\vd n}
  (-1)^{v(\lambda)}I_{\lambda/\alpha}^2=
  \sum_{\alpha/\mu\vd n}
  (-1)^{v(\mu)}I_{\alpha/\mu}^2
  $$
  if $n$ is even, and
  $$
  \sum_{\lambda/\alpha\vd n}
  (-1)^{v(\lambda)}I_{\lambda/\alpha}^2=
  \sum_{\alpha/\mu\vd n-1}
  (-1)^{v(\mu)}I_{\alpha/\mu}^2
  -\sum_{\alpha/\mu\vd n}
  (-1)^{v(\mu)}I_{\alpha/\mu}^2
  $$
  if $n$ is odd.
\end{theo}
\noindent
Figure~\ref{fig:bigexample} gives an example.
Observe that if $\alpha=\emptyset$ and $n\ge2$ the theorem reduces to
Theorem~\ref{th:jonasastrid2}.

\section{The proof of the main theorem\label{sec:mainproof}}
For a skew shape $\lambda/\mu$, let
$$\rsgn{\lambda/\mu}:=(-1)^{\sum_{(r,c)\in\lambda/\mu}(r-1)}.$$
For convenience, let $\rsgn T:=\rsgn{\sh T}$ for a skew tableau $T$.
Observe that for an ordinary shape $\lambda$ we have
$\rsgn\lambda=(-1)^{v(\lambda)}$.

For the sake of bookkeeping we will make two minor adjustments to
the skew insertion algorithm that do not affect the resulting
tableaux:
\begin{itemize}
\item
  Instead of starting with an empty Q-tableau, we
  start with the tableau $U$ after multiplying all entries by $\varepsilon$.
  Here $\varepsilon$ is a very small positive number.
\item
  During an internal
  insertion a new cell with an integer $b$ is added to the Q-tableau
  according to the usual rules. New additional rule: At the same time
  we remove the entry $b\varepsilon$ from the Q-tableau.
\end{itemize}

Consider the (adjusted) skew insertion algorithm starting with
P-tableau $P_0=T$ and Q-tableau $Q_0=U\varepsilon$.
After $\ell$ insertions (external or internal) we have obtained the
tableaux $P_\ell$ and $Q_\ell$. The following two lemmas state what
happens when we make the next insertion.
\begin{lemma}\label{lm:external}
  Let $(P_{\ell+1},Q_{\ell+1})$ be the resulting tableaux after
  {\em external} insertion of the number $a_1$ into $(P_\ell,Q_\ell)$.
  Then
  $$\frac{\sgn{P_{\ell+1}}}{\sgn{P_\ell}}
  =\frac{\sgn{Q_{\ell+1}}}{\sgn{Q_\ell}}
  \frac{\rsgn{Q_{\ell+1}}}{\rsgn{Q_\ell}}(-1)^{\sharp Q_\ell}(-1)^m,$$
  where $m$ is the number of entries in $P_\ell$ that are less than $a_1$.
\end{lemma}
\begin{proof}
  We insert the number $a_1$ which pops a number $a_2$ at $(1,c_1)$
  which pops a number $a_3$ at $(2,c_2)$ and so on. Finally the number
  $a_r$ fills a new cell $(r,c_r)$,
  see Figure~\ref{fig:extinsert}.
  \begin{figure}[t]
    \begin{center}
      \setlength{\unitlength}{5mm}
      \begin{picture}(14,12)(0,0)
        \put(0,0){\resizebox{70mm}{!}{\includegraphics{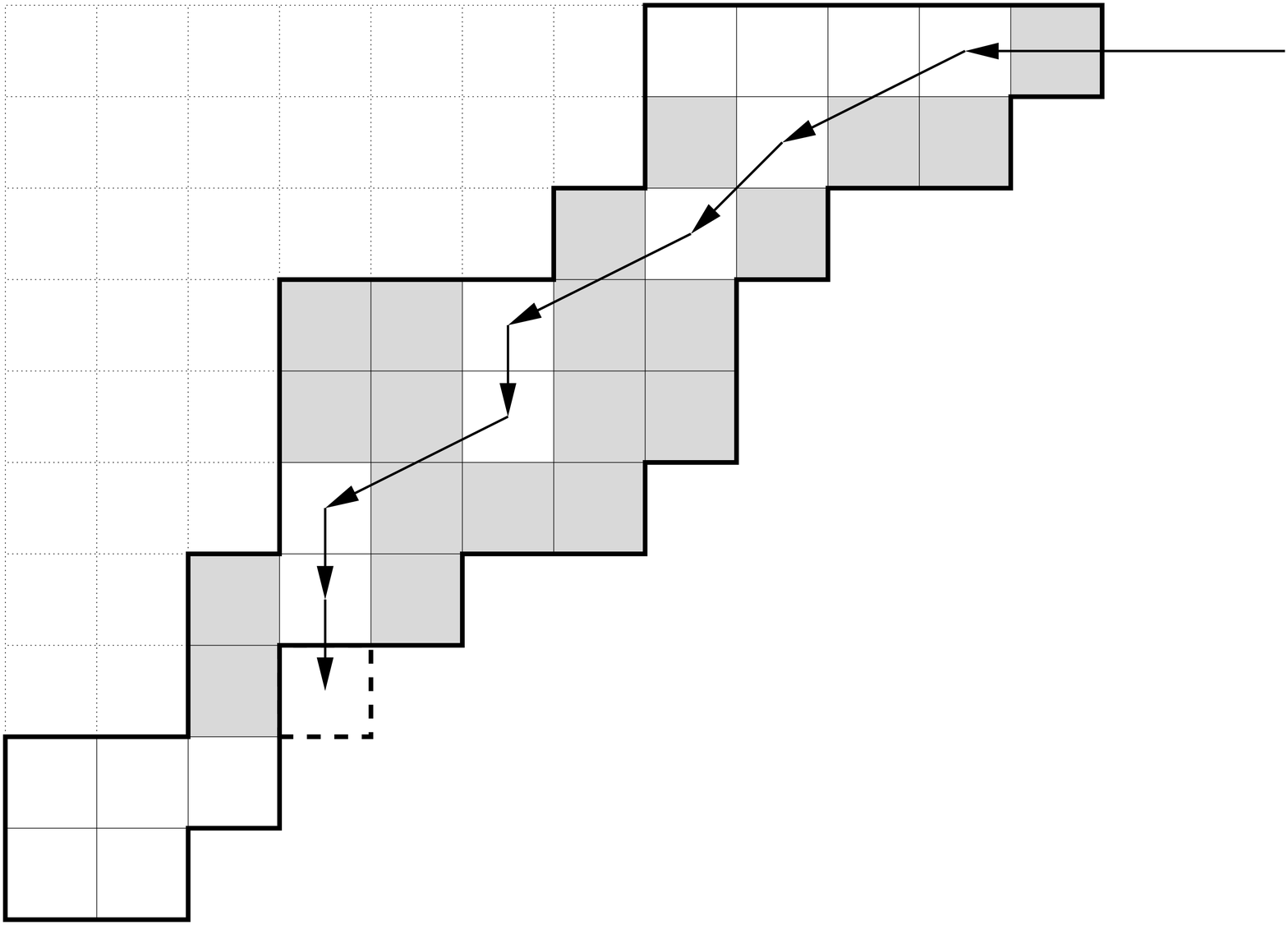}}}
        \put(14.4,9.4){$a_1$}
      \end{picture}
      \caption{External insertion of $a_1$. The shaded cells are
        counted by the sum
        $\sum_{i=2}^r(\beta_{i-1}-c_{i-1}+c_i-1-\gamma_i)$ in the proof.}
      \label{fig:extinsert}
    \end{center}
  \end{figure}
  
  For $2\le i\le r$, the relocation
  of $a_i$ multiplies the sign of the P-tableau by
  $(-1)^{\beta_{i-1}-c_{i-1}+c_i-1-\gamma_i}$,
  where $\sh{P_\ell}=\sh{Q_\ell}=\beta/\gamma$. Summation yields
  $$\sum_{i=2}^r(\beta_{i-1}-c_{i-1}+c_i-1-\gamma_i)
  =-(c_1-\gamma_1+r-2)+\sum_{i=1}^r(\beta_i-\gamma_i)$$
  since $\beta_r=c_r-1$.
  The placing of $a_1$ in the first row multiplies the sign of
  the P-tableau by $(-1)^{m-(c_1-1-\gamma_1)}$ where
  $m$ is the number of entries in $P_\ell$ that are
  less than $a_1$. We get
  $$\frac{\sgn{P_{\ell+1}}}{\sgn{P_\ell}}
  =(-1)^{m+1-r+\sum_{i=1}^r(\beta_i-\gamma_i)}.
  $$
  Obviously
  $$\frac{\invsgn{Q_{\ell+1}}}{\invsgn{Q_\ell}}
  =(-1)^{\sum_{i=1}^r(\beta_i-\gamma_i)}$$
  and
  $$\frac{\rsgn Q_{\ell+1}}{\rsgn Q_\ell}=(-1)^{r-1}.$$
  Since $\sgn R\invsgn R=(-1)^{\binom{\sharp R}{2}}$ for any
  tableau $R$, we have
  $$\frac{\invsgn Q_{\ell+1}}{\invsgn Q_\ell}
  =\frac{\sgn Q_{\ell+1}}{\sgn Q_\ell}(-1)^{\sharp Q_\ell}.$$
  Combining the equations above proves the lemma.
\end{proof}

\begin{lemma}\label{lm:internal}
  Let $(P_{\ell+1},Q_{\ell+1})$ be the resulting tableaux after
  {\em internal} insertion of the entry $a_1$ at $(r,c_0)$
  into $(P_\ell,Q_\ell)$. Then
  $$
  \frac{\sgn{P_{\ell+1}}}{\sgn{P_\ell}}
  =\frac{\sgn{Q_{\ell+1}}}{\sgn{Q_\ell}}
  \frac{\rsgn{Q_{\ell+1}}}{\rsgn{Q_\ell}}(-1)^{\sharp Q_\ell}.$$
\end{lemma}
\begin{proof}
  During an internal insertion the entry
  $a_1$ at $(r,c_0)$ pops a number $a_2$ at
  $(r+1,c_1)$ which pops a number $a_3$ at $(r+2,c_2)$
  and so on. Finally the number $a_k$ fills a new cell
  $(r+k,c_k)$, see Figure~\ref{fig:intinsert}.
  \begin{figure}[t]
    \begin{center}
      \setlength{\unitlength}{5mm}
      \begin{picture}(12,10)(0,0)
        \put(0,0){\resizebox{60mm}{!}{\includegraphics{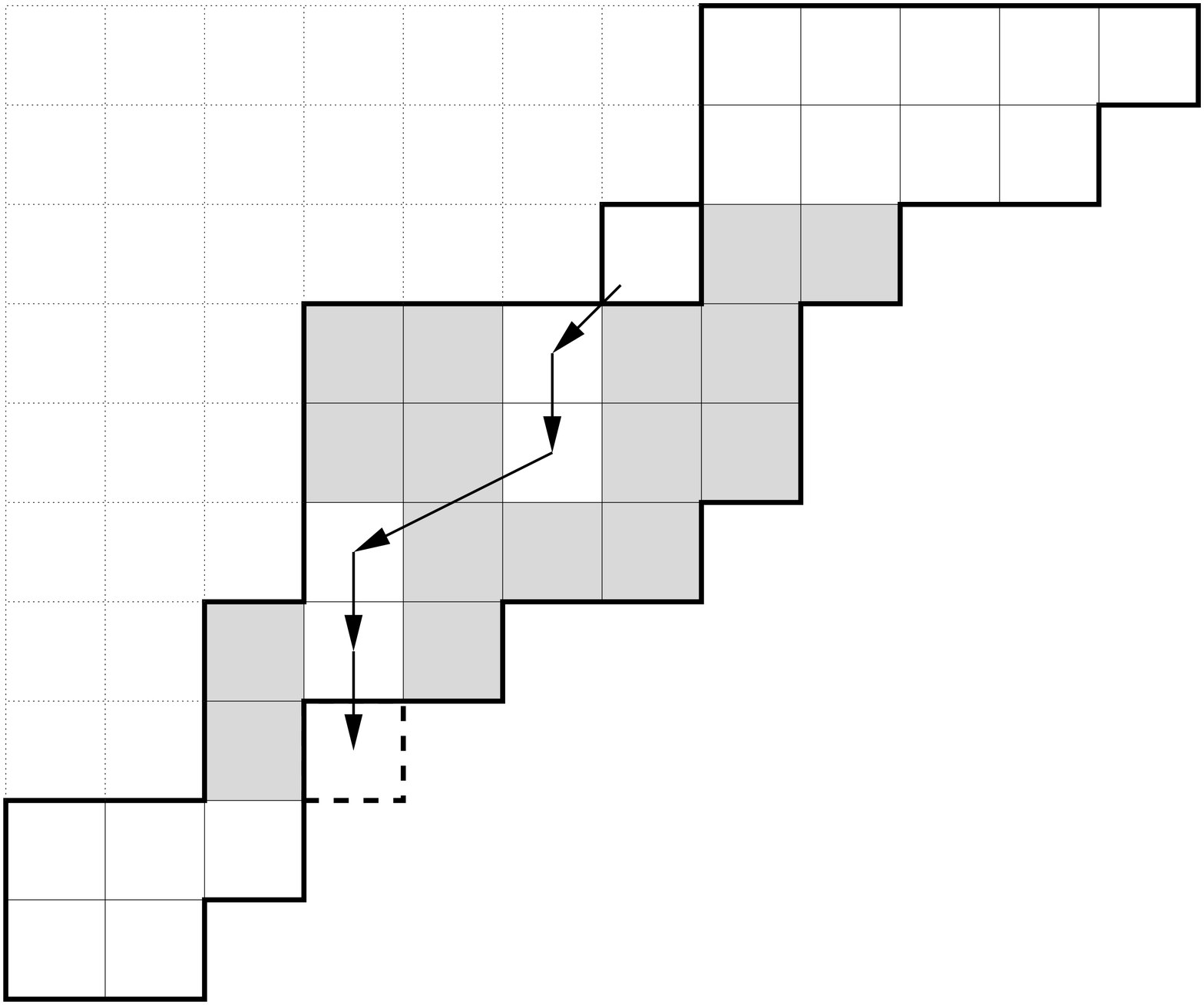}}}
        \put(6.2,7.3){$a_1$}
      \end{picture}
      \caption{Internal insertion starting with $a_1$. The shaded cells are
        counted by
        $\sum_{i=1}^k(\beta_{r+i-1}-c_{i-1}+c_{i}-1-\gamma_{r+i})$
        in the proof.}
      \label{fig:intinsert}
    \end{center}
  \end{figure}
  
  For $1\le i\le k$, the relocation of $a_i$ multiplies
  the sign of the P-tableau by
  $(-1)^{\beta_{r+i-1}-c_{i-1}+c_{i}-1-\gamma_{r+i}}$,
  where $\sh{P_\ell}=\sh{Q_\ell}=\beta/\gamma$. Summation yields
  $$\sum_{i=1}^k(\beta_{r+i-1}-c_{i-1}+c_{i}-1-\gamma_{r+i})
  =-k+\sum_{j=r}^{r+k}(\beta_j-\gamma_j)$$
  since $\beta_{r+k}=c_k-1$ and $\gamma_r=c_0-1$.
  
  What happens to the Q-tableau?
  According to our adjustments of the algorithm
  the entry $b\varepsilon$ at $(r,c_0)$ is removed and the entry $b$
  is added at the new cell at $(r+k,c_k)$. Observe that
  $b\varepsilon$ is the smallest element in $Q_\ell$;
  this is the very reason
  why we are making an internal insertion from its cell $(r,c_0)$.
  Also note that $b$ is the largest entry in $Q_{\ell+1}$. The transformation
  from $Q_\ell$ to $Q_{\ell+1}$ can be thought of as consisting of
  two steps: First we replace the entry $b\varepsilon$ by
  $b$, thereby changing the sign of the tableau by a factor
  $(-1)^{\sharp Q_\ell-1}$. Then we move the $b$ to the new cell at
  $(r+k,c_k)$, thereby changing the sign of the tableau by a factor
  $$(-1)^{-1+\sum_{j=r}^{r+k}(\beta_j-\gamma_j)}.$$
  Now, after observing that
  $$\frac{\rsgn{Q_{\ell+1}}}{\rsgn{Q_\ell}}=\frac{(-1)^{r+k}}{(-1)^r}
  =(-1)^k,$$
  the lemma follows.
\end{proof}
\vskip15mm
Now we are ready to prove our main theorem.
\begin{proof}[Proof of Theorem~\ref{th:main}.]
  From Lemma~\ref{lm:external} and~\ref{lm:internal}
  we deduce by induction that
  \begin{equation}\label{eq:intermed}
    \frac{\sgn P}{\sgn T}=\frac{\sgn Q}{\sgn U}\frac{\rsgn Q}{\rsgn U}
    (-1)^{\sum_{\ell=0}^{n-1}\sharp Q_\ell}(-1)^{\sum m}
  \end{equation}
  where $n=\sharp P$ and the last sum $\sum m$ is taken over
  all external insertions.

  Let $t_1<t_2<\cdots<t_g$ and
  $u_1<u_2<\cdots<u_g$ be the entries of $T$ and
  $U$, and write $\pi=\pair{i_1i_2\cdots i_h}{j_1j_2\cdots j_h}$.
  Let $\pi'$ be the permutation you get (in single-row notation)
  by preceding $\check\pi$ with
  the elements of $T$ decreasingly ordered,
  i.e.,~$\pi'=t_gt_{g-1}\cdots t_1j_1j_2\cdots j_h$.
  It is easy to see that
  the sum $\sum m$ equals the number of non-inversions of $\pi'$,
  i.e~pairs $i<j$ such that $\pi'(i)<\pi'(j)$. This means
  that $(-1)^{\sum m}=\invsgn{\pi'}$.

  What is the relationship between
  $\invsgn{\pi'}$ and $\sgn{\bar\pi}$?

  Let us go from $\pi'$ to $\bar\pi$ by a sequence
  of moves. Start with
  $$\pi'=t_gt_{g-1}\cdots t_1j_1j_2\cdots j_h.$$
  Move the first entry $t_g$ to position $u_g$:
  $$\underbrace{t_{g-1}t_{g-2}\cdots t_1\,j_\ast\cdots j_\ast\,t_g}
  _{\mbox{$u_g$ entries}}\,j_\ast\cdots j_\ast$$
  (Here the symbolic indices $\ast$ should be replaced by
  the sequence $1,2,\ldots,h$.)
  Next, move the entry $t_{g-1}$ to position $u_{g-1}$:
  $$\underbrace{t_{g-2}t_{g-3}\cdots t_1\,j_\ast\cdots j_\ast\,t_{g-1}}
  _{\mbox{$u_{g-1}$ entries}}\,j_\ast\cdots j_\ast\,t_g\,
  j_\ast\cdots j_\ast$$
  Continue until
  all elements of $T$ are moved. The resulting permutation is
  $\bar\pi$. After analysing what the moves do to the sign of the permutation,
  we obtain 
  $$\sgn{\bar\pi}=(-1)^K\sgn\pi'$$
  where
  $$K=\sum_{i=1}^g(u_i-1).$$
  Note also that
  $$\invsgn\pi'=(-1)^{\binom{n}{2}}\sgn\pi'.$$

  Now look at
  $$\sum_{\ell=0}^{n-1}\sharp Q_\ell.$$
  If we define $K_\ell:=\sharp\{b\in U\ :\ \ell<b\}$
  we can write $\sharp Q_\ell=\ell+K_\ell$. Summation yields
  $$\sum_{\ell=0}^{n-1}\sharp Q_\ell
  =\sum_{\ell=0}^{n-1}(\ell+K_\ell)=\binom{n}{2}+K+\sharp U.$$
  
  Now we are ready to update~(\ref{eq:intermed}):
  $$\frac{\sgn P}{\sgn T}=\frac{\sgn Q}{\sgn U}\frac{\rsgn Q}{\rsgn U}
  \sgn{\bar\pi}(-1)^{\sharp U}.$$
  There remains only some cleaning-up. Observe that
  $$\frac{\rsgn Q}{\rsgn U}
  =\frac{\rsgn{\lambda/\alpha}}{\rsgn{\alpha/\mu}}
  =\rsgn\lambda\rsgn\mu=(-1)^{v(\lambda)}(-1)^{v(\mu)}$$
  and $\sharp U=|\alpha|-|\mu|$.
  This yields the result
  $$(-1)^{v(\lambda)}\sgn P\sgn Q=(-1)^{|\alpha|}(-1)^{v(\mu)+|\mu|}
  \sgn T\sgn U\sgn{\bar\pi}.$$
\end{proof}

\section{The proof of Theorem~\ref{th:inout}\label{sec:inoutproof}}
In Theorem~\ref{th:ss} we have adopted the original notation from
Sagan and Stanley \cite{saganstanley}. However, for some applications
(and among them the forthcoming proof of Theorem~\ref{th:inout})
it is inconvenient to work with {\em partial} tableaux. For that matter
we now present a simple bijection that will allow us to work with
{\em standard} tableaux only.
\begin{lemma}\label{lm:seq}
  Let $n$ be a fixed positive integer and $\alpha$ and $\mu$ fixed
  partitions.
  Then there is a bijection $(\pi,T,U)
  \leftrightarrow(\tilde\pi,\tilde I,\tilde T,\tilde U)$ between
  \begin{itemize}
  \item
    triples $(\pi,T,U)$
    such that $\pi\in\PS_n$, $T,U\in\PT(\alpha/\mu)$
    and $\check{\pi}\cupdot T=\hat{\pi}\cupdot U=[n]$, and
  \item
    quadruples $(\tilde\pi,\tilde I,\tilde T,\tilde U)$ such that
    $\tilde\pi\in S_n$,
    $\tilde T,\tilde U\in\ST(\alpha/\mu)$ and $\tilde I\subseteq[n]$
    is the index set of an increasing subsequence of
    $\tilde\pi$ of length $|\alpha/\mu|$.
  \end{itemize}
  This bijection has the following properties:
  \begin{align*}
    \tilde\pi & =\bar\pi, \\
    \sgn{\tilde T} & =\sgn T, \\
    \sgn{\tilde U} & =\sgn U.
  \end{align*}
\end{lemma}
\begin{proof}
  Given a quadruple $(\tilde\pi,\tilde I,\tilde T,\tilde U)$,
  let the triple $(\pi,T,U)$ be given by the following procedure:
  Write $\tilde\pi$ in biword notation and remove the vertical pairs
  corresponding to the increasing subsequence $\tilde I$. The resulting
  partial permutation is $\pi$. Order the elements in $\tilde I$
  increasingly: $i_1<i_2<\cdots<i_k$. Now, for $1\le j\le k$,
  replace the entry $j$ in $\tilde U$ by $i_j$ and
  replace the entry $j$ in $\tilde T$ by $\tilde\pi(i_j)$. This results
  in $U$ and $T$ respectively. It is easy to see that this is indeed a
  bijection with the claimed properties.
\end{proof}

Now we are ready to prove Theorem~\ref{th:inout}.
\begin{proof}[Proof of Theorem~\ref{th:inout}.]
  Sum the equation of Theorem~\ref{th:main} over
  the whole domain of the skew RS-correspondence according to
  Theorem~\ref{th:ss} in view of Lemma~\ref{lm:seq}:
  \begin{gather*}
    \sum_{\lambda/\alpha\vd n}\;\sum_{P,Q\,\in\,\ST(\lambda/\alpha)}
    (-1)^{v(\lambda)}\sgn P\sgn Q=\\
    \sum_{k=0}^n\;\sum_{\alpha/\mu\vd k}\;
    \sum_{T,U\in\ST(\alpha/\mu)}\,\,
    \sum_{1\le i_1<\cdots<i_k\le n}\!\!\!\!
    \sum_{\sarray{
        \pi\in S_n \\
        \pi(i_1)<\cdots<\pi(i_k)
      }}
    \!\!\!\!\!\!\!\!\!\!(-1)^{|\alpha|+v(\mu)+|\mu|}
    \sgn T\sgn U\sgn\pi
  \end{gather*}
  Let $\LHS$ and $\RHS$ denote the left-hand side and
  the right-hand side of the equation above.
  The left-hand side trivially equals
  $$\LHS=\sum_{\lambda/\alpha\vd n}
  (-1)^{v(\lambda)}I_{\lambda/\alpha}^2.
  $$
  The right-hand side is trickier.
  Fix $1\le i_1<i_2<\cdots<i_k\le n$ and consider the sum
  $$S:=\sum_{\sarray{
      \pi\in S_n \\
      \pi(i_1)<\cdots<\pi(i_k)
    }}
  \!\!\!\!\!\!\sgn\pi.
  $$
  \begin{itemize}
  \item
    If $k=n$ clearly $S=1$.
  \item
    If $k\le n-2$ there are at least two integers
    $1\le a<b\le n$ not contained in the sequence
    $i_1<i_2<\cdots<i_k$. The sign-reversing involution
    $\pi\mapsto \pi\cdot(a,b)$ (here $(a,b)$ is
    the permutation that switches $a$ and $b$) shows that
    $S=0$.
  \item
    Suppose $k=n-1$ and let $a$ be the only integer in $[n]$
    not contained in the sequence $i_1<i_2<\cdots<i_k$.
    We are free to choose $\pi(a)$ from $[n]$, but as soon as
    $\pi(a)$ is chosen, the rest of $\pi$ must be the unique
    increasing sequence consisting of $[n]\setminus\pi(a)$ 
    if $\pi$ should contribute to $S$.
    The sign of $\pi$ then becomes $(-1)^{\pi(a)-a}$ so
    $$S=\sum_{i=1}^n(-1)^{i-a}=\left\{
      \begin{array}{ll}
        0 & \mbox{if $n$ is even,} \\
        (-1)^{a-1} & \mbox{if $n$ is odd.}
      \end{array}\right..$$
  \end{itemize}
  In the case where $n$ is odd and $k=n-1$, the double sum
  $$\sum_{1\le i_1<\cdots<i_k\le n}\!\!\!
  \sum_{\sarray{
      \pi\in S_n \\
      \pi(i_1)<\cdots<\pi(i_k)
    }}
  \!\!\!\!\!\!\sgn\pi
  =\sum_{a=1}^n(-1)^{a-1}=1.$$
  
  In summary we have showed
  $$\sum_{1\le i_1<\cdots<i_k\le n}\!\!\!
  \sum_{\sarray{
      \pi\in S_n \\
      \pi(i_1)<\cdots<\pi(i_k)
    }}
  \!\!\!\!\!\!\sgn\pi
  =\left\{\begin{array}{rl}
      1 & \mbox{if $k=n$,}\\
      1 & \mbox{if $k=n-1$ and $n$ is odd,}\\
      0 & \mbox{if $k=n-1$ and $n$ is even,}\\
      0 & \mbox{if $k\le n-2$.}
    \end{array}\right.$$
  
  If $n$ is even we finally obtain
  $$\RHS=(-1)^{|\alpha|}\;\sum_{\alpha/\mu\vd n}\;
  (-1)^{v(\mu)+|\mu|}
  \sum_{T,U\in\ST(\alpha/\mu)}
  \sgn T\sgn U
  =\sum_{\alpha/\mu\vd n}
  (-1)^{v(\mu)}I_{\alpha/\mu}^2
  $$
  since $(-1)^{|\alpha|+|\mu|}=(-1)^{|\alpha|-|\mu|}=(-1)^n=1$.

  Analogously, if $n$ is odd we get
  \begin{align*}
    \RHS & =(-1)^{|\alpha|}\;\sum_{\alpha/\mu\vd n}\;
    (-1)^{v(\mu)+(|\alpha|-n)}
    \sum_{T,U\in\ST(\alpha/\mu)}
    \sgn T\sgn U\\
    & \quad\mbox{}+(-1)^{|\alpha|}\;\sum_{\alpha/\mu\vd n-1}\;
    (-1)^{v(\mu)+(|\alpha|-(n-1))}
    \sum_{T,U\in\ST(\alpha/\mu)}
    \sgn T\sgn U\\
    & =\sum_{\alpha/\mu\vd n-1}
    (-1)^{v(\mu)}I_{\alpha/\mu}^2
    -\sum_{\alpha/\mu\vd n}
    (-1)^{v(\mu)}I_{\alpha/\mu}^2.
  \end{align*}
\end{proof}

\section{Specialisations of Theorem~\ref{th:inout}\label{sec:specialisations}}
Apart from the special case $\alpha=\emptyset$,
Theorem~\ref{th:inout} offers a couple of other nice
specialisations if we choose the parameters $\alpha$ and $n$ properly.
First we obtain a surprising formula for the square of
the sign-imbalance of any ordinary shape:
\begin{coro}
  Let $\alpha$ be a fixed $n$-shape. Then
  $$I_\alpha^2
  =\sum_{\lambda/\alpha\vd n}
  (-1)^{v(\lambda)}I_{\lambda/\alpha}^2
  =\sum_{\lambda/\alpha\vd n+1}
  (-1)^{v(\lambda)}I_{\lambda/\alpha}^2$$
  if $n$ is even, and
  $$I_\alpha^2
  =\sum_{\lambda/\alpha\vd n-1}
  (-1)^{v(\lambda)}I_{\lambda/\alpha}^2
  $$
  if $n$ is odd.
\end{coro}
\begin{proof}
  First suppose $n$ is even.
  Theorem~\ref{th:inout} yields
  $$
  \sum_{\lambda/\alpha\vd n}
  (-1)^{v(\lambda)}I_{\lambda/\alpha}^2=
  \sum_{\alpha/\mu\vd n}
  (-1)^{v(\mu)}I_{\alpha/\mu}^2.
  $$
  The right-hand side consists of only one term, namely
  $(-1)^{v(\emptyset)}I_{\alpha/\emptyset}^2=I_\alpha^2$.
  From Theorem~\ref{th:inout} we also get
  $$
  \sum_{\lambda/\alpha\vd n+1}
  (-1)^{v(\lambda)}I_{\lambda/\alpha}^2=
  \sum_{\alpha/\mu\vd n}
  (-1)^{v(\mu)}I_{\alpha/\mu}^2
  -\sum_{\alpha/\mu\vd n+1}
  (-1)^{v(\mu)}I_{\alpha/\mu}^2.
  $$
  The second term of the right-hand side vanishes and the first term
  is $I_\alpha^2$ as before.

  Now suppose $n$ is odd. Then Theorem~\ref{th:inout} yields
  $$
  \sum_{\lambda/\alpha\vd n-1}
  (-1)^{v(\lambda)}I_{\lambda/\alpha}^2=
  \sum_{\alpha/\mu\vd n-1}
  (-1)^{v(\mu)}I_{\alpha/\mu}^2
  $$
  The right-hand side consists of only one term, namely
  $(-1)^{v((1))}I_{\alpha/(1)}^2$ which equals $I_\alpha^2$
  since in an ordinary tableau the 1 is always located at $(1,1)$.
\end{proof}

Next we present another generalization of Theorem~\ref{th:jonasastrid2}.
\begin{coro}
  Let $\alpha$ be a fixed $n$-shape. Then
  $$\sum_{\lambda/\alpha\vd m}
  (-1)^{v(\lambda)}I_{\lambda/\alpha}^2=0$$
  for any integer $m\ge n+2$ if $n$ is even, and for
  any integer $m\ge n$ if $n$ is odd.
\end{coro}
\begin{proof}
  If $m$ is even Theorem~\ref{th:inout} yields
  $$
  \sum_{\lambda/\alpha\vd m}
  (-1)^{v(\lambda)}I_{\lambda/\alpha}^2=
  \sum_{\alpha/\mu\vd m}
  (-1)^{v(\mu)}I_{\alpha/\mu}^2.
  $$
  The right-hand side vanishes since $m>|\alpha|$.

  If $m$ is odd Theorem~\ref{th:inout} yields
  $$
  \sum_{\lambda/\alpha\vd m}
  (-1)^{v(\lambda)}I_{\lambda/\alpha}^2=
  \sum_{\alpha/\mu\vd m-1}
  (-1)^{v(\mu)}I_{\alpha/\mu}^2
  -\sum_{\alpha/\mu\vd m}
  (-1)^{v(\mu)}I_{\alpha/\mu}^2.
  $$
  If $m\ge n+2$ the right-hand side vanishes simply because
  $m-1>|\alpha|$. Otherwise $n$ is odd and the only remaining case
  is $m=n$. But then the right-hand side becomes
  $I_{\alpha/(1)}^2-I_\alpha^2=0$.
\end{proof}

\section{Future research\label{sec:future}}
For an ordinary shape $\lambda$, let $h(\lambda)$ be the
number of disjoint horizontal dominoes that fit in $\lambda$
and let $d(\lambda)$ be the number of disjoint
$2\times2$-squares (fourlings) that fit in $\lambda$.

In~\cite{jonas} the following theorem, conjectured
by Stanley \cite{stanley},
was proved (the (a)-part was independently proved by T.~Lam~\cite{lam}):
\begin{theo}[Stanley; Lam; Sj\"ostrand; 2003]\label{th:stanleyab}
\mbox{\\}
\begin{enumerate}
\item[(a)] For every $n\geq0$
$$\sum_{\lambda\vdash n}q^{v(\lambda)}t^{d(\lambda)}x^{h(\lambda)}
I_\lambda=(q+x)^{\lfloor n/2\rfloor}.$$
\item[(b)] If $n\not\equiv1\pmod4$
$$\sum_{\lambda\vdash n}(-1)^{v(\lambda)}t^{d(\lambda)}I_\lambda^2=0.$$
\end{enumerate}
\end{theo}
The (b)-part is a strengthening of Theorem~\ref{th:jonasastrid2}
and one might wonder if there is a similar strengthening
of Theorem~\ref{th:inout} for skew shapes.

The (a)-part
is about signed sums of sign-imbalances without taking
the square. From an RS-correspondence perspective
it is unnatural not to take the square of the sign-imbalance
since the P- and Q-tableaux come in pairs. In fact it
might be argued that non-squared sign-imbalances are unnatural
in all cases, because their sign is dependent on the actual
labelling of the poset, i.e.,~it is important that
we read the tableau as a book. However, part (a)
in the theorem is still true (and there are even stronger theorems,
see~\cite{jonas}) and it can be proved by means of the
RS-correspondence as was done in~\cite{jonas}. This suggests
that the skew RS-algorithm could be a useful tool for
studying signed sums of non-squared sign-imbalances too.

As a tool for proving Theorem~\ref{th:stanleyab} the concept of
{\em chess tableaux} was introduced in~\cite{jonas}. A chess tableaux
is a standard Young tableau where odd entries are located at an
even Manhattan distance from the upper-left cell of the shape,
while even entries are located at odd distances.
This notion of course generalizes to skew tableaux
(in fact it generalizes to many
other posets) and since it proved so useful in the study of
sign-imbalance of ordinary shapes we think it will shed some light
on the skew shapes as well.

Another direction of research is to find analogues to
Theorem~\ref{th:main} for other variants of the RS-algorithm.
For instance, in~\cite[Theorem~5.1]{saganstanley}
Sagan and Stanley present a
generalization of their skew RS-correspondence where
the condition that $\sh U=\sh T$ and $\sh P=\sh Q$ is
relaxed. From that they are able to infer identities like
$$\sum_\sarray{
\lambda/\beta\vd n \\
\lambda/\alpha\vd m
}
f_{\lambda/\beta}f_{\lambda/\alpha}
=\sum_{k\ge0}\binom{n}{k}\binom{m}{k}k!
\sum_\sarray{
\alpha/\mu\vd n-k \\
\beta/\mu\vd m-k
}
f_{\alpha/\mu}f_{\beta/\mu}$$
where $f_{\lambda/\mu}=\sharp\ST(\lambda/\mu)$. This correspondence
may give interesting formulas for sums of products of
sign-imbalances as well.

\section*{Acknowledgements}
I want to thank Richard Stanley and Frank Sottile for
personal communication.

\end{document}